\documentclass[12pt]{amsart}

\usepackage{amsmath,amssymb}

\newtheorem{theorem}{Theorem}[section]
\newtheorem{lemma}[theorem]{Lemma}

\theoremstyle{definition}
\newtheorem{proposition}[theorem]{Proposition}

\newtheorem{corollary}[theorem]{Corollary}
\newtheorem{example}[theorem]{Example}

\theoremstyle{remark}
\newtheorem{remark}[theorem]{Remark}

\numberwithin{equation}{section}

\let\ord=\preccurlyeq
 
\def\det{\mathop{\mathbf{det}}}
\def\tr{\mathop{\mathbf{tr}}}

\begin{document}

\title[Hilbert's metric on symmetric cones]{Application of Hilbert's projective metric on symmetric cones }

\author{Khalid KOUFANY}
\address{
Institut Elie Cartan,
UMR 7502 (UHP-CNRS-INRIA),
Universit{\'e} H. Poincar{\'e}, B.P. 239,
F-54506 Vand{\oe}uvre-l{\`e}s-Nancy Cedex,
France
}
\email{khalid.koufany@iecn.u-nancy.fr}





\begin{abstract}
Let $\Omega$ be a symmetric cone. In this note, we introduce Hilbert's projective
metric on $\Omega$ in terms of Jordan algebras and we apply it to prove
that given a linear transformation $g$ such
that $g(\Omega)\subset \Omega$ and a real number $p$, $|p|>1$, then there
exists a unique element $x\in\Omega$ satisfying $g(x)=x^p$.
\end{abstract}

\maketitle

\section{Introduction}
If $A$ is a positive linear transformation on $\mathbb{R}^n$, then the Perron Theorem says that there exists $x_0\in\mathbb{R}^n_+\setminus\{0\}$ such that for all $x\in\mathbb{R}^n_+$, $A^nx$ converges in direction to $x_0$, i.e. $\frac{A^nx}{\|A^nx\|}\to\frac{x_0}{\|x_0\|}$.
It has been shown by Birkhoff \cite{Birkhoff} that the Perron Theorem can be considered  as a special case of the Banach contraction theorem. Birkhoff's approach is to consider a projective metric in a cone of positive elements in a Banach space. This metric was first introduced by Hilbert \cite{Hilbert} when studying  hyperbolic geometry.  We begin with the definition of this metric in a general setting.

Let $V$ be is a real Banach space and $C$ be a closed convex pointed  cone, where pointed means that $C\cap -C=\{0\}$. Denote the corresponding interior by $\Omega$. The relation $\preccurlyeq$ is defined on $V$ by saying that $x\ord y$ if and only if $y-x\in C$. Hence $(V,\ord)$ is a partially ordered linear space and it is Archimedean, that is, if $ny\ord y$ for all $n\in\mathbb{N}^*$, then $y\ord 0$.\\
For $x\in V$ and $y\in \Omega$ we let
\begin{equation*}
M(x,y):=\inf\{\lambda \;|\; x\ord \lambda y\},
\end{equation*}
and
\begin{equation*}
m(x,y):=\sup\{\mu \;|\; \mu y\ord x\}.
\end{equation*}
Hilbert's projective metric is defined on $\Omega$ by
\begin{equation}\label{hilbert_metric}
d(x,y)=\log\frac{M(x,y)}{m(x,y)}.
\end{equation}
In the case of $\mathbb{R}^n_+$, Hilbert's projective metric is
$d(x,y)=\log\frac{\mathrm{max}{\frac{x_i}{y_i}}}{\mathrm{min}{\frac{x_i}{y_i}}}$
where $x=(x_1,\ldots,x_n)$ and $y=(y_1,\ldots,y_n)$ are two vectors of
$\mathbb{R}^n_+$.

The Hilbert projective metric may be applied to variety of problems
involving positive matrices and positive integral operators. For example
one can use it to solve some Volterra equations. It is also particularly
useful in proving the existence of the fixed point for positive
operators defined in a Banach space. In this way, it has been shown by
Bushell \cite{Bushell3} that Hilbert projective metric may be applied to
prove that, if $T$ a real nonsingular $r\times r$ matrix, then there exists
a unique real positive definite symmetric $r\times r$ matrix $A$ such that 
\begin{equation}\label{equa_bushell}
T'AT=A^2.
\end{equation}
Notice that if $T$ is neither symmetric nor orthogonal the existence and
the uniqueness of $A$  is not an elementary problem, even if $r=3$. 

In this note we will formulate the Hilbert projective metric on symmetric cones in a way most convenient for our purpose using Jordan algebra theory
and extend Bushell's Theorem to this class of convex cones.

\section{Preparatory definitions}\label{preliminaries}

We summarize here a number of basic definitions and results concerning symmetric cones. For more detailed discussion, see Faraut and Kor{\'a}nyi \cite{Faraut-Koranyi}.

Let $\Omega$ be an open convex cone in a Euclidean vector space $V$ of dimension $n$. Let $G(\Omega)$ be the group of linear automorphisms  of $\Omega$
$$G(\Omega)=\{g\in GL(V) \;|\; g(\Omega)=\Omega\}.$$
Then $\Omega$ is said to be {\it homogeneous} if $G(\Omega)$ acts on it transitively. If $\overline{\Omega}$ is pointed, then $\Omega$ is said to be {\it symmetric} if it is homogeneous and self-dual.\\
A {\it Euclidean Jordan algebra} is a Euclidean vector space $V$ equipped with a bilinear product such that
\begin{align*}
  xy&=yx\\
  x(x^2y)&=x^2(xy)\\
  (xy|z)&=(y|xz).
\end{align*}
It is shown that the interior $\Omega$ of the set of squares in $V$ is a symmetric cone, and every symmetric cone is given in this way.\\
We define the (left) multiplication $L$ by $L(x)y=xy$ and the so-called {\it quadratic representation} $P$ by $P(x)=2L^2(x)-L(x^2)$. For any $x\in V$, the endomorphisms $L(x)$ and $P(x)$ are self-adjoint.
\begin{example}\label{example1}
The space $\mathrm{Sym}(r,\mathbb{R})$ of real symmetric
$r\times r$
matrices is a Jordan algebra for the
product $x\circ y=\frac{1}{2}(xy+yx)$ and it is Euclidean for the scalar
product $(x|y)=\mathrm{Tr}(xy)$. The associated symmetric cone  is the set
$\Omega_\mathrm{Sym}$ of positive definite symmetric $r\times r$ matrices and  $G(\Omega_\mathrm{Sym})$ is the linear group $\mathrm{GL}(r,\mathbb{R})$.
In this case, $P(x)y=xyx$.
\end{example}
The Euclidean Jordan algebra has and identity element denoted by $e$. Let $K(\Omega):=\{g\in G(\Omega) \;|\; g(e)=e\}=$, then $K(\Omega)$ is a maximal compact subgroup of $G(\Omega)$. 
Let $r$ be the rank of $V$. A {\it Jordan frame} $\{c_1,\ldots,c_r\}$ of $V$ is a complete system of non-zero orthogonal primitive idempotents :
$$
\begin{aligned}
{}&c_i^2=c_i, \; c_i\;\mathrm{indecomposable},\\
{}&c_ic_j=0\;\mathrm{if}\; i\not=j,\\
{}&c_1+\ldots+c_r=e.
\end{aligned}
$$
We denote by $\mathcal{J}(V)$ the set of all primitive idempotents of $V$.
Suppose $V$ is simple. In other words,  there is no non-trivial ideal in $V$. Then each element $x$ in $V$ can be written as
\begin{equation}\label{polar}
x=k(\sum_{j=1}^r\lambda_jc_j),\; k\in K(\Omega)_\circ,\; \lambda_j\in\mathbb{R},
\end{equation}
where $K(\Omega)_\circ$ is the identity component of $K(\Omega)$.
The {\it determinant} is defined by $\det(x)=\prod_{j=1}^r\lambda_j$ and the {\it trace form} by $\tr(x)=\sum_{j=1}^r\lambda_j$. The real numbers $\lambda_1,\ldots,\lambda_r$ are the {\it eigenvalues} of $x$
and the {\it spectral norm} of $x$ is then defined by
\begin{equation}\label{spectral_norm}
|x|=\sup_i|\lambda_i|.
\end{equation}
In the case of $\mathrm{Sym}(r,\mathbb{R})$, (\ref{polar}) corresponds to
the polar decomposition (diagonalization) of symmetric matrices. The
functions $\det$ and
$\tr$ are respectively the usual determinant $\mathrm{Det}$ and trace $\mathrm{Tr}$ of matrices.\\

\section{Hilbert's metric on symmetric cones}
If we consider the cone $\Omega_\mathrm{Sym}$ of real symmetric positive
definite $r\times r$ matrices, then one can easily express the Hilbert
projective metric (\ref{hilbert_metric}) in terms of eigenvalues of
elements of $\Omega_\mathrm{Sym}$. Indeed, if $A$ and $B$ are in $\Omega_\mathrm{Sym}$, then 
$$M(A,B):=\inf\{\lambda \;|\; \lambda B-A\ord 0\}=\max_{\|x\|=1}\frac{(Ax|x)}{(Bx|x)},$$
and 
$$m(A,B):=\sup\{\lambda \;|\; \lambda B-A\ord 0\}=\min_{\|x\|=1}\frac{(Ax|x)}{(Bx|x)},$$ which are respectively the
greatest and the least eigenvalue of $B^{-1}A$. Observe that eigenvalues of
the matrix $B^{-1}A$ are the same of
matrix $B^{-\frac{1}{2}}AB^{-\frac{1}{2}}=P(B^{-\frac{1}{2}})A$ (see Example 
\ref{example1} for this notation).

More generally, for symmetric cones, Hilbert's projective metric can be also formulated in terms of extremal eigenvalues : let $x$ and $y$ be in $\Omega$ and let $\lambda_M(x,y)>0$ and $\lambda_m(x,y)>0$ denote the greatest and the least eigenvalue of the element $P(y^{-\frac{1}{2}})x\in \Omega$. Then one can prove the following (see  \cite[Thoerem 4.2]{Koufany2})
\begin{proposition}\label{min_max} We have
\begin{equation}
\lambda_M(x,y)=\max_{c\in\mathcal{J}(V)}\frac{(x|c)}{(y|c)},
\end{equation}
and
\begin{equation}
\lambda_m(x,y)=\min_{c\in\mathcal{J}(V)}\frac{(x|c)}{(y|c)}.
\end{equation}
\end{proposition}

From this characterization one can easily prove the following
\begin{lemma}\label{lm_Brauer}
If $x,\, y\in\Omega$ and $\alpha, \,\beta\geq 0$, then
\begin{enumerate}
\item $\lambda_M(\alpha x+\beta y,y)=\alpha\lambda_M(x,y)+\beta$,
\item $\lambda_m(\alpha x+\beta y,y)=\alpha\lambda_m(x,y)+\beta$,
\item $\lambda_M(x,y){\lambda_m(y,x)}=1$.
\end{enumerate}
\end{lemma}

\begin{proposition}\label{def_hil}
If $x,\, y\in\Omega$, then the Hilbert metric of $x$ and $y$ is given by
\begin{equation}
d(x,y)=\log\frac{\lambda_M(x,y)}{\lambda_m(x,y)}=\log\lambda_M(x,y)\lambda_M(y,x).
\end{equation}
\end{proposition}
\begin{proof}
According to Proposition \ref{min_max}, we have
$$\begin{array}{ll}
M(x,y)&=\inf\{\lambda \;|\; x\ord\lambda y\}\\
      &=\inf\{\lambda \;|\; \lambda y-x\in\overline{\Omega}\}\\
      &=\max\{\frac{(x|c)}{(y|c)} \;|\; c\in\mathcal{J}(V)\}\\
      &=\lambda_M(x,y),
\end{array}
$$
and
$$\begin{array}{ll}
m(x,y)&=\sup\{\mu \;|\; \mu y\ord\lambda x\}\\
      &=\sup\{\mu \;|\; x-\mu y\in\overline{\Omega}\}\\
      &=\min\{\frac{(x|c)}{(y|c)} \;|\; c\in\mathcal{J}(V)\}\\
      &=\lambda_m(x,y).
\end{array}
$$
Applying Lemma \ref{lm_Brauer}, we obtain easily the second equality of the proposition.
\end{proof}

\begin{remark}{\rm
It follows from Lemma \ref{lm_Brauer} and Proposition \ref{def_hil} that $d$ is constant on rays,
 \begin{equation}\label{const_rays}
d(\lambda x,\mu y)=d(x,y)\;\;\mathrm{for}\; \lambda,\,\mu>0.
 \end{equation}
In particular, if $\mu>0$, then the map $\xi_\mu : x\mapsto \mu x$ is an isomery of $d$.
}
\end{remark}

\section{Contractions of Hilbert's projective metric}
Let $x=\sum_{j=1}^r\lambda_j c_j$ be an element of $\Omega$. If
$p\in\mathbb{R}$, we write $x^p$ for $\sum_{j=0}^r\lambda_j^pc_j$ and we
state  the following contraction result :
\begin{proposition} Let $p\in\mathbb{R}$ such that $|p|\leq 1$. Then the map
\begin{equation}\label{power_p}
\omega_{p} : \Omega\to \Omega,\; x \mapsto x^p 
\end{equation}
is a contraction of $d$, that is, for any $x,\, y\in\Omega$,
\begin{equation}\label{contraction}
d(x^p,y^p)\leq |p|\,d(x,y).
\end{equation}
\end{proposition}
\begin{proof}
Let $x, \,y\in\Omega$. Then $\lambda_M(x^{-1},y^{-1})=1/\lambda_m(x,y)$ and $\lambda_m(x^{-1},y^{-1})=1/\lambda_M(x,y)$. Hence the map 
\begin{equation}\label{imath}
\imath : \Omega \to \Omega : x \mapsto x^{-1}
\end{equation}
is an isometry of $d$.  
It suffices then to prove (\ref{contraction}) for $0\leq p\leq 1$.
We know that $x\ord \lambda_M(x,y)y$, then by the Loewner-Heinz inequality,
see \cite[Corollary 9]{Lim1}, we have $x^p\ord
\lambda_M(x,y)^py^p$. Therefore $\lambda_M(x^p,y^p)\leq
\lambda_M(x,y)^p$. Using the same arguments we prove that $\lambda_m(x^p,y^p)\geq \lambda_m(x,y)^p$. Thus $d(x^p,y^p)\leq p\, d(x,y)$.
\end{proof}

\section{Completeness}
\begin{proposition}
$(\Omega,d)$ is a pseudo-metric space. In other words, for any $x,\, y,\, z\in\Omega$, the following holds :
\begin{enumerate}
\item[(a)] $d(x,y)\geq 0$
\item[(b)] $d(x,y)=d(y,x)$
\item[(c)] $d(x,z)\leq d(x,y)+d(y,z)$
\item[(d)] $d(x,y)=0\Leftrightarrow \exists \lambda>0 : x=\lambda y$
\end{enumerate}
\end{proposition}
\begin{proof}
 Let $x, y$ and $z$ be in $\Omega$. The first propriety of the proposition is obvious, since $\lambda_M(x,y)\geq \lambda_m(x,y)$.\\
(b) From Lemma \ref{lm_Brauer} we obtain
$$\begin{array}{ll}
d(x,y)&=\log\frac{\lambda_M(x,y)}{\lambda_m(x,y)}\\
      &=\log\lambda_M(x,y)\lambda_M(y,x)\\
      &=d(y,x).
\end{array}$$
(c) Let $c$ be a primitive idempotent of $V$, then according to Proposition \ref{min_max}, $0<\frac{(x|c)}{(y|c)}\leq \lambda_M(x,y)$ and $0<\frac{(y|c)}{(z|c)}\leq \lambda_M(y,z)$. Hence $\frac{(x|c)}{(z|c)}=\frac{(x|c)}{(y|c)} \frac{(y|c)}{(z|c)}\leq \lambda_M(x,y)\lambda_M(y,z)$ and 
\begin{equation}\label{ineq_tri_0}
\lambda_M(x,z)\leq\lambda_M(x,y)\lambda_M(y,z).
\end{equation} 
Similarly we prove that 
\begin{equation}\label{ineq_tri}
\lambda_m(x,z)\geq \lambda_m(x,y)\lambda_m(y,z). 
\end{equation} 
It follows then from (\ref{ineq_tri_0}) and (\ref{ineq_tri}) that $d(x,z)\leq
d(x,y)+d(y,z)$.\\
(d) If $d(x,y)=0$, then $\lambda_M(,x,y)=\lambda_m(x,y):=\lambda$ and all
the eigenvalues of $P(y^{-\frac{1}{2}})x$ are equal to $\lambda$. Therefore, $P(y^{-\frac{1}{2}})x=\lambda e$ and $x=\lambda P(y^{\frac{1}{2}})e=\lambda y$. Conversely, if $x=\lambda y$ where $\lambda >0$, then $P(y^{-\frac{1}{2}})x=\lambda e$ and $\lambda_M(x,y)=\lambda_m(x,y)=\lambda$.
\end{proof}

\medskip

Let $S(V)$ be the unite sphere in $V$, $S(V)=\{x\in V \;|\; |x|=1\}$ with respect to the spectral norm introduced in (\ref{spectral_norm}).

\begin{lemma}\label{lm_cle1} Let $x$ and $y$ be in $\Omega\cap S(V)$. Then
\begin{enumerate} 
\item[(a)]  $|x-y|\leq e^{d(x,y)}-1$.
\item[(b)]  If  $|x-y|<\lambda_m(y)$ then  $|x-y|\geq \lambda_m(y)\tanh\bigl\{\frac{1}{2}d(x,y)\bigr\}$, where $\lambda_m(y)$ is the least eigenvalue of $y$.
\end{enumerate}
\end{lemma}
\begin{proof}
Let $x, y\in \Omega\cap S(V)$ and let $c$ be a primitive idempotent of $V$, then by Proposition (\ref{min_max}) we have
\begin{equation}\label{3_ineg}
\lambda_m(x,y)(y|c)\leq (x|c)\leq \lambda_M(x,y)(y|c).
\end{equation}
But $|x|=1$ and $|y|=1$, then there exists $c_1, c_2\in \mathcal{J}(V)$ such that $(x|c_1)=1$ and $(y|c_2)=1$, thus 
\begin{equation}\label{3_ineg_bis}
\lambda_m(x,y)\leq 1 \leq \lambda_M(x,y).
\end{equation}
Let $c_0$ be a primitive idempotent of $V$ such that $|x-y|=(x-y|c_0)$,
then using (\ref{3_ineg_bis}) we have
$$\begin{array}{ll}
|x-y|&=(x|c_0)-(y|c_0)\\
     &=\bigl\{\frac{(x|c_0)}{(y|c_0)}-1\bigr\} (y|c_0)\\ 
     &\leq \frac{(x|c_0)}{(y|c_0)}-1 \;\;\; \text{since}\;\; |y|=1\\
     &\leq \lambda_M(x,y)-1\\
     &\leq \lambda_M(x,y)-\lambda_m(x,y) \\
     &\leq \bigl\{\lambda_M(x,y)-\lambda_m(x,y)\bigr\}\frac{1}{\lambda_m(x,y)}\\
     &=e^{d(x,y)}-1.
\end{array}
$$
Moreover
\begin{eqnarray}\label{ineq1_lemma}
\lambda_M(x,y)&&=\max_{c\in\mathcal{J}(V)}\frac{(x|c)}{(y|c)} \nonumber\\
              &&=\max_{c\in\mathcal{J}(V)}\{\frac{(x-y|c)}{(y|c)}+1\}\nonumber\\
              &&\leq |x-y|/\lambda_m(y)+1.
\end{eqnarray}

Now, if $|x-y|<\lambda_m(y)$, then $1-\frac{|x-y|}{\lambda_m(y)}>0$. Therefor, from (\ref{ineq1_lemma})  and Lemma \ref{lm_Brauer}, we
have $\bigl(1-\frac{|x-y|}{\lambda_m(y)}\bigr)y\ord x$. Hence
\begin{equation}\label{ineq2_lemma}
\lambda_m(x,y)\geq 1-\frac{|x-y|}{\lambda_m(y)}.
\end{equation}
Finally, it follows from (\ref{ineq1_lemma}) and (\ref{ineq2_lemma}) that $|x-y|\geq \lambda_m(y)\tanh\{\frac{1}{2}d(x,y)\}$.
\end{proof}

\begin{proposition}\label{complete}
$(\Omega\cap S(V),d)$ is complete metric space.
\end{proposition}
\begin{proof}
 It is clear that if $x,\, y\in(\Omega\cap S(V),d)$ such that $d(x,y)=0$,
 then $x=y$. Thus $(\Omega\cap S(V),d)$ is a metric space. \\
Let $(x_k)_k$ be a Cauchy sequence in $(\Omega\cap S(V),d)$. Then from
Lemma \ref{lm_cle1}, $(x_k)_k$ is a Cauchy sequence in $(V, |\cdot|)$ and then converges to an element $x\in \overline{\Omega}\cap S(V)$. 
In order to prove that $x$ is an element of $\Omega\cap S(V)$, we observe,
using (\ref{ineq_tri}) that $\lambda_m(x_i)=\lambda_m(x_i,e)\geq \lambda_m(x_i,x_j)\lambda_m(x_j,e)=\lambda_m(x_j)\lambda_M(x_i,x_j)e^{-d(x_i,x_j)}$.
Therefor, from (\ref{3_ineg_bis}), we obtain $\lambda_m(x_i)\geq
\lambda_m(x_j)e^{-d(x_i,x_j)}$. Since $(x_k)_k$ converges in the spectral norm
to $x$, $(\lambda_m(x_k))_k$ converges to $\lambda_m(x)$. It follows that for a fixed
large $j$, $\lambda_m(x)=\lim_{i\to\infty}\lambda_m(x_i)\geq
\alpha\lambda_m(x_j)>0$, where $\alpha>0$. Hence $x\in \Omega\cap
S(V)$. 
Finally, it is easy to prove using the second proposition of Lemma \ref{lm_cle1} that $(x_k)$ converges to $x$ in $(\Omega\cap S(V),d)$.
\end{proof}

\section{Application}
If $V$ is the Jordan algebra $\mathrm{Sym}(r,\mathbb{R})$, then $\Omega$ is the symmetric cone of positive definite matrices, and $G(\Omega)$ coincides with the linear group $\mathrm{GL}(r,\mathbb{R})$. Recall that in this case $G(\Omega)$ acts on $\Omega$  by  : $t\cdot a=t'at$.\\
In \cite{Bushell3}, Bushell prove that for any $t\in \mathrm{GL}(r,\mathbb{R})$ and any integer $k\geq 1$, there exists a unique real positive definite matrix $a$ such that $t'at=a^{2^k}$. \\
For a general symmetric cone we have the following 

\begin{theorem}
Let $g\in G(\Omega)$ and $p\in\mathbb{R}$ such that $|p|> 1$. Then there exists a unique element $a$ in $\Omega$ such that 
$g(a)=a^p$. 
\end{theorem}
\begin{proof}
Since $\imath : x\mapsto x^{-1}$ is an isometry of $d$, see (\ref{imath}),
we can assume that $p>1$. For $x\in \Omega\cap S(V)$ we put
$F(x)=\frac{1}{|g(x)^{\frac{1}{p}}|}g(x)^{\frac{1}{p}}$. Then $F$ maps
$\Omega\cap S(V)$ into $\Omega\cap S(V)$. Recall that Hilbert's projective
metric is constant on rays, see (\ref{const_rays}), then we have for $x,
y\in\Omega\cap S(V)$, 
$$d(F(x),F(y))=d((\omega_{1/p}\circ g)(x),(\omega_{1.p}\circ g)(y),$$
where $\omega_{1/p}$ is the power map $x\mapsto x^{\frac{1}{p}}$. We have
already noticed in (\ref{power_p}) that $\omega_{1/p}$ is a contraction of
the metric $d$, hence
$$d(F(x),F(y))\leq \frac{1}{p} d(g(x),g(y)).$$
Moreover, the eigenvalues of $P(y^{-1/2})x$ are the unique solutions of the characteristic equation $\det(\lambda y-x)=0$. But, using
\cite[Proposition III.4.3]{Faraut-Koranyi}, 
$$\det(\lambda
g(y)-g(x))=\det(g(\lambda y-x))=\mathrm{Det}(g)^{n/r}\det(\lambda
y-x).$$ Thus  the eigenvalues $P(y)^{-1/2}x$ and $P(g(y))^{-1/2}g(x)$ are
the same and $g$ is also an isometry of $d$. Finally all this implies that
$F$ is a contraction of the Hilbert metric $d$,
$$d(F(x),F(y))\leq \frac{1}{p} d(x,y).$$
 Therefore by the Banach contraction mapping theorem  and Proposition
 \ref{complete} the map $F$ has a unique fixed point, let say $u$. Then $u^p=\frac{1}{|g(u)^{1/p}|^p}g(u)$ and the element $a:=|g(u)^{1/p}|^pu\in\Omega$ is the unique solution of the equation $g(x)=x^p$.
\end{proof}

\begin{corollary}
Let $h\in G(\Omega)$ and $p\in\mathbb{R}$ such that $|p|>1$. Then there exists a unique element $a$ of $\Omega$ such that 
$h(a^p)=a$. 
\end{corollary}

\noindent{\bf Final remark.} The Hilbert original definition of the projective metric involved the logarithm of the cross-ratio of for points in the interior of a convex cone in $\mathbb{R}^n$. It would  be interesting to define the Hilbert projective metric on symmetric cones using the generalized cross-ratio introduced in \cite{Braun}, see also \cite{Koufany2}.

\bibliographystyle{amsplain}

\begin{thebibliography}{10}

\bibitem{Birkhoff}
 G. Birkhoff,
\textit{Extensions of Jentzsch's theorem.}
Trans. Amer. Math. Soc.,
\textbf{85} (1957), 219-227. 
{\tt MR19:296a }

\bibitem{Braun}
H. Braun,
\textit{Doppelverh{\"a}ltnisse in Jordan-Algebren.}
Abh. Math. Sem. Univ. Hamburg
\textbf{32} (1968), 25--51

\bibitem{Bushell2}
 P. J. Bushell,
\textit{Hilbert's metric and positive contraction mappings in a Banach space.}
Arch. Rational Mech. Anal.,
\textbf{52} (1973), 330-338.
{\tt MR 49:1247}

\bibitem{Bushell3}
 P. J. Bushell,
\textit{On solutions of the Matrix Equation $T'AT=A^2$.}
Linear Algebra and Appl.,
\textbf{8} (1974), 465--469. 
{\tt MR 51:551}

\bibitem{Faraut-Koranyi}
 J. Faraut and A. Kor{\'a}nyi, 
\textit{Analysis on Symmetric Cones.}
Oxford Mathematical Monographs, Clarendon Press, Oxford, 
(1994).
{\tt MR 98g:17031}

\bibitem{Hilbert} 
D. Hilbert,
\textit{{\"U}ber die gerade Linie als k{\"u}rzeste Verbindung zweier Punkte.}
Math. Ann.,
\textbf{46} (1895), 91--96.

\bibitem{Koufany2}
 K. Koufany,
\textit{Contractions of angles in symmetric cones.}
Publ. RIMS, Kyoto Univ. 
\textbf{38} (2002), n$^\mathrm{o}$ 2.



\bibitem{Lim1}
Y. Lim,
\textit{Applications of geometric means on symmetric cones.}
Math. Ann.,
\textbf{319} (2001), 457--468. 
{\tt MR 1819877 }


\end{thebibliography}

\end{document}